\documentclass[10pt]{amsart}
\usepackage[applemac]{inputenc}
\usepackage{amsaddr}

\usepackage[colorlinks=true,linkcolor=red,citecolor=blue,urlcolor=blue,hypertexnames=false]{hyperref}
\usepackage{bookmark}
\usepackage{amsthm,thmtools,amssymb,amsmath,amscd}

\usepackage{fancyhdr}
\usepackage{esint}
\usepackage{enumerate}

\usepackage{pictexwd,dcpic}

\swapnumbers
\declaretheorem[name=Theorem,numberwithin=section]{thm}
\declaretheorem[name=Remark,style=remark,sibling=thm]{rem}
\declaretheorem[name=Lemma,sibling=thm]{lemma}

\declaretheorem[name=Definition,style=definition,sibling=thm]{defn}
\declaretheorem[name=Corollary,sibling=thm]{cor}

\declaretheorem[name=Example,style=remark,sibling=thm]{example}

\numberwithin{equation}{section}

\usepackage{cleveref}
\crefname{lemma}{Lemma}{Lemmata}
\crefname{prop}{Proposition}{Propositions}
\crefname{thm}{Theorem}{Theorems}
\crefname{cor}{Corollary}{Corollaries}
\crefname{defn}{Definition}{Definitions}
\crefname{example}{Example}{Examples}
\crefname{rem}{Remark}{Remarks}
\crefname{assum}{Assumption}{Assumptions}
\crefname{nota}{Notation}{Notation}

\theoremstyle{remark}

\usepackage{autonum}

\newcommand{\ti}{\tilde}

\newcommand{\bs}{\backslash}
\newcommand{\cn}{\colon}
\newcommand{\sub}{\subset}

\newcommand{\bbR}{\mathbb{R}}

\newcommand{\bbS}{\mathbb{S}}

\newcommand{\8}{\infty}

\newcommand{\al}{\alpha}

\newcommand{\e}{\epsilon}
\newcommand{\ka}{\kappa}
\newcommand{\la}{\lambda}

\newcommand{\s}{\sigma}

\newcommand{\Om}{\Omega}

\newcommand{\G}{\Gamma}

\newcommand{\La}{\Lambda}

\newcommand{\cL}{\mathcal{L}}

\newcommand{\cD}{\mathcal{D}}

\newcommand{\cF}{\mathcal{F}}

\newcommand{\del}{\partial}
\newcommand{\n}{\nabla}

\newcommand{\fa}{\forall}
\newcommand{\rt}{\sqrt}

\newcommand{\ip}[2]{\left\langle #1,#2 \right\rangle}
\newcommand{\fr}[2]{\frac{#1}{#2}}
\newcommand{\x}{\times}

\DeclareMathOperator{\dive}{div}
\DeclareMathOperator{\id}{id}

\DeclareMathOperator{\osc}{osc}

\DeclareMathOperator{\dist}{dist}

\DeclareMathOperator{\tr}{tr}
\DeclareMathOperator{\Rm}{Rm}

\DeclareMathOperator{\grad}{grad}

\DeclareMathOperator{\vol}{vol}
\DeclareMathOperator{\Area}{Area}
\DeclareMathOperator{\EV}{EV}

\newcommand{\pf}[1]{\begin{proof}{\parskip\baselineskip{ #1}} \end{proof}}
\newcommand{\eq}[1]{\begin{equation}\begin{alignedat}{2} #1 \end{alignedat}\end{equation}}

\newcommand{\br}[1]{\left(#1\right)}
\newcommand{\abs}[1]{\lvert #1\rvert}
\newcommand{\enum}[1]{\begin{enumerate}[(i)] #1 \end{enumerate}}

\newcommand{\ra}{\rightarrow}

\newcommand{\mt}{\mapsto}

\newcommand{\mc}{\mathcal}

\newcommand{\mrm}{\mathrm}

\newcommand{\hp}{\hphantom}
\newcommand{\q}{\quad}

\begin{document}
\title[Curvature flows at Matrix 2020]{Extrinsic curvature flows and applications}

\author[J. Scheuer]{Julian Scheuer}
\address{Columbia University New York/Universität Freiburg}
\email{\href{mailto:julian.scheuer@math.uni-freiburg.de}{julian.scheuer@math.uni-freiburg.de}}

\date{\today}
\thanks{These lectures were held at the {\it{Early Career Researchers Workshop on Geometric Analysis and PDEs}} at the Matrix Institute in Creswick, Australia. I would like to thank the Matrix Institute and the organizers Paul Bryan, Jiakun Liu, Mariel Saéz and Haotian Wu for the opportunity to give these lectures.}

\begin{abstract}
These notes arose from a mini lecture series the author gave at the Early Career Researchers Workshop on Geometric Analysis and PDEs, held in January 2020 at the Matrix institute of the University of Melbourne.
We discussed some classical aspects of expanding curvature flows and obtained first applications.
 In these notes we will give a detailed account on what was covered during the lectures.
\end{abstract}

%
%

\maketitle

\section{Introduction}

\subsection*{Expanding curvature flows}
This is an introduction to the theory of (expanding) extrinsic curvature flows, i.e. normal variations of hypersurfaces the speed of which are determined by the principal curvatures at each point. The flowing hypersurfaces are parametrized by a time-dependent family of embeddings
\eq{x\cn [0,T)\x \bbS^{n}\ra \bbR^{n+1}}
which satisfies
\eq{\label{intro:flow}\dot{x}=\fr{1}{f(\ka_{1},\dots,\ka_{n})}\nu,}
where 
\eq{\ka_{1}\leq \dots\leq \ka_{n}}
are the principal curvatures at $x$, $\nu$ is the outward pointing unit normal and a dot denotes the partial time derivative.

Under a monotonicity assumption on $f$, this flow is a weakly parabolic system and we present proofs of the classical results due to Claus Gerhardt \cite{Gerhardt:/1990} and John Urbas \cite{Urbas:/1990}: Under certain assumptions on $f$ and the initial embedding $x_{0}$ this flows exists for all times and after exponential blowdown converges to a round sphere. Furthermore we show that this flow can be used to prove so-called {\it{Alexandrov-Fenchel inequalities}}, which are inequalities between certain curvature functionals of a hypersurface. The approach is due to Pengfei Guan and Junfang Li \cite{GuanLi:08/2009}. Classical examples are the isoperimetric inequality and the Minkowski inequality
\eq{\int_{M}H\geq c_{n}\abs{M}^{\fr{n-1}{n}},}
which holds if $M$ is mean-convex ($H>0$) and starshaped. Here $\abs{M}$ is the surface area of $M$. Equality holds precisely on every geodesic sphere. An appropriate rescaling of the flow \eqref{intro:flow} has nice monotonicity properties which, together with the convergence result, can be used to prove the inequalities. The approach we take slightly differs from the original works \cite{Gerhardt:/1990,GuanLi:08/2009}. Namely we use that the normal component of the rescaled flow actually moves by 
\eq{\label{LC}\dot x=\br{\fr{1}{f}-\fr{u}{n}}\nu,}
where $u$ is the support function of the hypersurface. A priori estimates for \eqref{intro:flow} are directly deduced along this rescaling, which makes the estimates a little easier compared to Gerhardt's original arguments \cite{Gerhardt:/1990}.
One interesting aspect of this particular rescaling is that \eqref{LC} belongs to the class of so-called {\it{locally constrained curvature flows}}. The mean curvature type flow of this class,
\eq{\dot x=(n-uH)\nu,}
was invented by Pengfei Guan and Junfang Li in \cite{GuanLi:/2015} as a natural flow to prove the isoperimetric inequality in space forms: It preserves the enclosed volume and decreases surface area. A variety of such flows have appeared since then and they have been useful to obtain new geometric inequalities, cf. \cite{GuanLiWang:/2019,HuLiWei:02/2020,ScheuerWangXia:11/2018,ScheuerXia:11/2019,WangXia:10/2019}.

\subsection*{Outline}
These notes are structured as follows. First we present some background on the curvature function $f$. It is known that the ordered principal curvatures are continuous in time, but if they have higher multiplicity they are in general not smooth. Hence at first sight the operator in \eqref{intro:flow} seems to lack regularity. However, this issue can be worked around by considering the function
\eq{F(A)=f\circ \EV(A),}
where $A$ is the Weingarten (or shape-) operator of the embedding and $\EV$ the eigenvalue map. Interestingly, even though $\EV$ is not smooth, if $f$ is smooth and symmetric, $F$ will be a smooth and natural map on the space of vector space endomorphisms. To people working with fully nonlinear curvature operators this is well known. We will give the precise setup to make this approach rigorous and state some important relations between derivatives of $f$ and $F$, but skip most of the proofs in these notes. The material is taken from \cite{Scheuer:06/2018}.

Afterwards we first fix some notation and conventions about hypersurface geometry and deduce the evolution equations for various geometric quantities. After these general considerations, we actually start with the a priori estimates for the inverse curvature flows and prove their convergence.
We conclude by presenting the application to Alexandrov-Fenchel inequalities.

Up to some hard results from general parabolic PDE theory, i.e. short-time existence of fully nonlinear equations, Krylov-Safonov- and Schauder theory, the exposition should be mostly self-contained. However, on some occasions we will skip proofs for elementary statements.

\section{Curvature functions}\label{sec:CF}
We quickly introduce the algebra of curvature functions using a new approach from \cite{Scheuer:06/2018}.
Along a variation
\eq{\dot{x}=-f\nu}
the function $f$ is supposed to be a function of the principal curvatures of the flow hypersurfaces $M_{t}=x(t,M)$. As we deal with geometric flows, $f$ has to be invariant under coordinate changes and thus we require it to be symmetric under all permutations. Hence we may assume the $\ka_{i}$ to be ordered,
\eq{\ka_{1}\leq \dots\leq \ka_{n}.} 
We assume that $f$ is smooth. Along the curvature flows considered later, we derive estimates for the curvature and hence we would like to deduce a parabolic equation which is satisfied by the $\ka_{i}$. However, those are in general not smooth functions, so we need to find another description of $f$, namely make it depend on the Weingarten operator $A$, the components of which are smooth.

This can be accomplished with the following idea: Suppose $\G\sub\bbR^{n}$ is an open and symmetric set and
\eq{f\in C^{\8}(\G)}
symmetric.
It is a classical result \cite{Glaeser:01/1963} that $f$ then is a function of the elementary symmetric polynomials
\eq{\label{ES}s_{m}(\ka):=\sum_{1\leq i_{1}<\dots<i_{m}\leq n}\prod_{j=1}^{m}\ka_{i_{j}},}
or also of the power sums
\eq{p_{m}(\ka)=\sum_{i=1}^{n}\ka_{i}^{m}.}
This means
\eq{f=\rho(s_{1},\dots,s_{n})=\psi(p_{1},\dots,p_{n})}
for some smooth functions $\rho$ and $\psi$.
The crucial point is, that for the power sums it is very easy to make the transition from the dependence on the eigenvalues $\ka_{i}$ to dependence on the operator. This is formalized as follows:

{\defn{Let $V$ be an $n$-dimensional real vector space and $\cD(V)\sub \cL(V)$ be the set of real diagonalizable endomorphisms. Then we denote by $\mrm{EV}$ the eigenvalue map, i.e.
\eq{\mrm{EV}\cn \cD(V)&\ra \bbR^{n}/\mc{P}_{n}\\
				A&\mt (\ka_{1},\dots,\ka_{n}),}
where $\ka_{1},\dots, \ka_{n}$ denote the eigenvalues of $A$ and $\mc{P}_{n}$ is the permutation group of $n$ elements.
}}

For the power sums there is a very obvious candidate to serve as a function defined on linear maps, namely
\eq{P_{k}(A)=\tr(A^{k}).}
Then there holds
\eq{P_{k}(A)=p_{k}(\mrm{EV}(A))\q\fa A\in \cD(V).}
Now we can just insert the $P_{k}$ into $\psi$, i.e. we  define
\eq{F=\psi(P_{1},\dots,P_{n}).}
Then $F\in C^{\8}(\Om)$ for some open set $\Om \sub \cL(V)$ and
\eq{F_{|\cD_{\G}(V)}=f\circ\EV_{|\cD_{\G}(V)},}
where $\cD_{\G}(V)$ is the set of those real diagonalizable linear maps with eigenvalues in $\G$. We obtain the following relations for the derivatives, see \cite{Scheuer:06/2018} for the details. Denote by $F'(A)$ the gradient of $F$, i.e. by the relation
\eq{dF(A)B=\tr(F'(A)\circ B).} If $A$ is real diagonalizable, then $F'(A)$ is real diagonalizable and if we denote by $F^{i}(A)$ its eigenvalues, then
\eq{F^{i}(A)=\fr{\del f}{\del\ka_{i}}(\ka),}
where $\ka=\EV(A)$. The second derivatives are related via
\eq{\label{2DerHk-A}d^{2}F(A)(\eta,\eta)=\sum_{i,j=1}^{n}\fr{\del^{2} f}{\del\ka_{i}\del\ka_{j}}\eta^{i}_{i}\eta_{j}^{j}+\sum_{i\neq j}^{n}\fr{\fr{\del f}{\del \ka_{i}}-\fr{\del f}{\del\ka_{j}}}{\ka_{i}-\ka_{j}}\eta^{i}_{j}\eta^{j}_{i},}
where $f$ is evaluated at the $n$-tuple $(\ka_{i})$ of corresponding eigenvalues. The latter quotient is also well defined in case $\ka_{i}=\ka_{j}$ for some $i\neq j$. Here $(\eta^{i}_{j})$ is a matrix representation of some $\eta\in \cL(V)$ with respect to a basis of eigenvectors of $A$.

Later we will require $F$ to have certain properties, which we collect in the following definition.

\begin{defn}
The function $F$ is called
\enum{
\item homogeneous of degree one, if $\G$ is a cone and 
\eq{F(\la A)=\la F(A)\q\fa \la>0~\fa A\in \cD_{\G}(V),}
\item strictly monotone, if
	\eq{\EV(F'(A))\in \G_{+}\q\fa A\in \cD_{\G}(V),}
\item concave, if
\eq{D^{2}F(A)(\eta,\eta)\leq 0}
for all $A$ and for all $\eta$ which are jointly self-adjoint with $A$.
}
Here $\G_{+}$ is the positive open cone on $\bbR^{n}$,
\eq{\G_{+}=\{\ka\in \bbR^{n}\cn \ka_{i}>0\q\fa 1\leq i\leq n\}.}
\end{defn}

\begin{example}
Important examples of functions $f$, such that $F$ has the properties in the above definition, are the quotients
\eq{q_{m}=\fr{s_{m}}{s_{m-1}}}
or the roots
\eq{\s_{m}=s_{m}^{\fr 1m}.}
In either case $f$ has the mentioned properties in the cone
\eq{\G_{m}=\{\ka\in \bbR^{n}\cn s_{k}>0\q \fa 1\leq k\leq m\},}
see for example \cite{HuiskenSinestrari:09/1999}.
Later we will use the quotients to deduce the Alexandrov-Fenchel inequalities.
\end{example}

\section{Some hypersurface geometry}
\subsection{Conventions on Riemannian geometry}
In this section we state the basic conventions concerning the elementary objects of Riemannian geometry.
Let $M$ be a smooth manifold of dimension $n$. For vector fields $X, Y$ which are also derivations of $C^{\8}(M)$, their Lie bracket is given by
\eq{[X,Y]=XY-YX}
and for an endomorphism field $A$ we denote by $\tr A\in C^{\8}(M)$ its trace.
Let $g$ be a Riemannian metric on $M$ with Levi-Civita connection $\n$.
The Riemannian curvature tensor is
\eq{\Rm(X,Y)Z=\n_{X}\n_{Y}Z-\n_{Y}\n_{X}Z-\n_{[X,Y]}Z,}
and we also use $\Rm$ to denote the associated $(0,4)$-tensor,
\eq{\Rm(X,Y,Z,W)=g(\Rm(X,Y)Z,W).}
The connection $\n$ induces covariant derivatives of tensor fields $T$ in the usual way via
\eq{&~\n T(X_{1},\dots,X_{l},Y^{1},\dots,Y^{k},X)\\
=&\br{\n_{X}T}(X_{1},\dots,X_{l},Y^{1},\dots,Y^{k})\\
		=&~X(T(X_{1},\dots,X_{l},Y^{1},\dots,Y^{k}))-T(\n_{X}X_{1},X_{2},\dots,X_{l},Y^{1},\dots,Y^{k})\\
		&-\ldots-T(X_{1},\dots,X_{l},Y^{1},\dots,\n_{X}Y^{k}).}
Let
\eq{x\cn M\ra \bbR^{n+1}}
be the smooth embedding of an $n$-dimensional manifold. The induced metric of $x(M)$ is given by the pullback of the ambient Euclidean metric $\ip{\cdot}{\cdot}$,
\eq{g=x^{*}\ip{\cdot}{\cdot}.}
 The second fundamental form $h$ of the embedding $x$ is given by the Gaussian formula
\eq{\label{Gauss}D_{x_{\ast}(X)}x_{\ast}(Y)=x_{\ast}(\n_{X}Y)- h(X,Y)\nu,}
where $D$ is the standard Euclidean connection.
The Weingarten operator is defined via
\eq{g(A(X),Y)=h(X,Y)}
and the Weingarten equation says that
\eq{\label{Weingarten}D_{x_{\ast}(X)}\nu=x_{\ast}(A(X)).}
Finally, we have the Gauss equation,
\eq{\label{GaussEq}\Rm(W,X,Y,Z)&=h(W,Z)h(X,Y)-h(W,Y)h(X,Z).}
	
\begin{rem}\label{Rem-coordinates}
We will simplify the notation by using the following shortcuts occasionally:
\begin{enumerate}[(i)]
\item We will often omit $x_{\ast}$, i.e. when we insert a tangent vector field $X$ into an ambient tensor field, we always understand $X$ to be given by its pushforward. 
\item When we deal with complicated evolution equations of tensors, we will occasionally use a local frame to express tensors with the help of their components, i.e. for a $(k,l)$-tensor field $T$, an expression like $T^{i_{1}\dots i_{k}}_{j_{1}\dots j_{l}}$ is understood to be
\eq{T^{i_{1}\dots i_{k}}_{j_{1}\dots j_{l}}=T(e_{j_{1}},\dots,e_{j_{l}},\e^{i_{1}},\dots \e^{i_{k}}),}
where $(e_{i})$ is a local frame and $(\e^{i})$ its dual coframe.
\item The coordinate expression for  the $m$-th covariant derivative of a $(k,l)$-tensor field $T$ is
\eq{\n^m T=\br{\n_{j_{l+m}\dots j_{l+1}}T^{i_1\dots i_k}_{j_1\dots j_l}},}
where subscripts to $\n$ represent the derivatives.
\end{enumerate}
\end{rem}

\subsection{Hypersurfaces in polar coordinates}
The punctured Euclidean space is isometric to
\eq{N=(0,\8)\x\bbS^{n}, \q\bar g=dr^{2}+r^{2}\s,}
where $\s$ is the round metric on $\bbS^{n}$ and $r=\abs{x}$. We will deal with closed starshaped hypersurfaces, i.e. those which can be written as graphs over the fibre $\bbS^{n}$. We collect some useful formulae here.

 Differentiating twice along $M$ and using the Gaussian formula \eqref{Gauss} gives
\eq{\label{ThetaDeriv}\fr 12\n^{2}\abs{x}^{2}=g-uh,}
where $u$ is the {\it{support function}}
\eq{\label{support}u=\ip{r\del_{r}}{\nu}=\ip{x}{\nu}.}

The flow hypersurfaces we consider are graphs over $\bbS^{n}$, so let us recall some standard formulae, which can be found in \cite[Sec.~1.5]{Gerhardt:/2006}. Let $M_{0}=x(M)\sub \bbR^{n+1}$ be a graph over $\bbS^{n}$,
\eq{M_{0}=\{(\rho(y),y)\cn y\in \bbS^{n}\}=\{(\rho(y(\xi)),y(\xi))\cn \xi\in M\}.}
Then the induced metric of $M_{0}$ is
\eq{g=d\rho \otimes d\rho+\rho^{2}\s.}
We choose the normal $\nu$ to satisfy
\eq{\ip{\nu}{\del_{r}}>0.}
Let
\eq{\bar h=\rho\s}
 be the second fundamental form of the embedded slice $\{r=\rho\},$
then the second fundamental form of $M_{0}$ can be expressed with the help of the graph function,
\eq{\label{graph-h}uh=-\rho\n^{2}\rho+\rho\bar{h}=-\rho\n^{2}\rho+g-d\rho\otimes d\rho,}
 which is an easy exercise using the Gaussian formula and the Christoffel-symbols in polar coordinates.
 Also note that the principal curvatures $\bar\ka$ of these slices are given by
 \eq{\bar\ka=\fr{1}{\rho}.}

\subsection*{Formulae for hypersurface variations}
As we consider time-dependent families of embedded hypersurfaces, we have to know how the previously discussed geometric quantities behave along variations with arbitrary speed,
\eq{\dot{x}=-\cF\nu,}
where $\nu$ is the same normal as the one in the Gaussian formula \eqref{Gauss}.
\begin{lemma}\label{Ev}
Let $T>0$, $M^{n}$ a smooth orientable manifold and
\eq{x\cn [0,T)\x M\ra \bbR^{n+1}}
be a normal variation with velocity $-\cF$ of a smooth hypersurface $M_{0}=x(0,M)$. Then the following evolution equations are satisfied.
\begin{enumerate}[(i)]
\item The induced metric $g$ satisfies
\eq{\label{Ev-g}\dot{g}=-2 \cF h.}
\item The normal vector field satisfies
\eq{\label{Ev-nu}\fr{D}{dt}\nu=\grad \cF,}
where $\tfrac{D}{dt}$ is the covariant time derivative along the curve $x(\cdot,\xi)$ for fixed $\xi\in M$.
\item The Weingarten operator evolves by
\eq{\label{Ev-W-1}\dot{A}=\n\grad \cF+ \cF A^{2}.}

\end{enumerate}
\end{lemma}

\pf{Let $X,Y$ be vector fields.

``\eqref{Ev-g}'': Due to the Weingarten equation \eqref{Weingarten} we have
 \eq{\dot g(X,Y)&=\ip{D_{\dot x}X}{Y}+\ip{X}{D_{\dot{x}}Y}=- \cF\ip{D_{X}\nu}{Y}- \cF\ip{X}{D_{Y}\nu}=-2 \cF h(X,Y).}
				
``\eqref{Ev-nu}'': We have
\eq{0=\fr{\del}{\del t}\ip{\nu}{\nu}=\ip{\fr{D}{dt}\nu}{\nu}}
and
\eq{\ip{\fr{D}{dt}\nu}{X}=-\ip{\nu}{D_{\dot{x}}X}=X\cF=\ip{\grad \cF}{X}.}

``\eqref{Ev-W-1}'': Differentiate the Weingarten equation \eqref{Weingarten} with respect to time. The left hand side gives
\eq{D_{\dot{x}}D_{X}\nu&=D_{X}D_{\dot{x}}\nu=\n_{X}\grad \cF- h(X,\grad \cF)\nu,}
where we have used \eqref{Ev-nu}. The right hand side gives
\eq{D_{\dot{x}}(A(X))&=D_{A(X)}\dot{x}+\dot{A}(X)=- h(X,\grad \cF)\nu- \cF A^{2}(X)+\dot{A}(X).}
Equate both sides to get the result.
}

\section{Classical inverse curvature flows}\label{ICF}

We prove the classical result of Claus Gerhardt \cite{Gerhardt:/1990} and John Urbas \cite{Urbas:/1990}, that the inverse curvature flow
\eq{\dot{x}=\fr{1}{F}\nu}
in the Euclidean space $\bbR^{n+1}$,
starting from starshaped and $F$-admissable\footnote{At every point the Weingarten operator is in the domain of definition} initial data converges to a round sphere after rescaling. Here is the result in detail.

\begin{thm}[\cite{Gerhardt:/1990,Urbas:/1990}]\label{thm:ICF}
Let $n\geq 2$ and $x_{0}\in C^{\8}(\bbS^{n},\bbR^{n+1})$ be the embedding of a starshaped $F$-admissable hypersurface, where $F\in C^{\8}(\G)\cap C^{0}(\bar\G)$ is a positive, strictly monotone, $1$-homogeneous and concave curvature function on a symmetric, open and convex cone $\G$ which contains $(1,\dots,1)$. Suppose that
\eq{F_{|\G}>0,\q F_{|\del\G}=0,\q F(1,\dots,1)=n.}  Then the parabolic Cauchy-problem
\eq{\label{ICF}\dot{x}&=\fr{1}{F}\nu\\
		x(0,\cdot)&=x_{0}}
has a unique solution $x\in C^{\8}([0,\8)\x\bbS^{n},\bbR^{n+1}).$
The rescaled hypersurfaces
\eq{\ti x(t,\cdot)=e^{-\fr tn}x(t,\cdot)}
converge smoothly to the embedding of a round sphere.
\end{thm}

We use an approach slightly different from the original papers, namely we work directly on the rescalings. Note that $\ti x$ will solve
\eq{\label{Resc}\dot{\ti x}=\fr{1}{F(e^{\fr tn}A)}\ti \nu-\fr{1}{n}\ti x.}
As the Weingarten operator scales reciprocally to the hypersurfaces,
\eq{\ti A=e^{\fr tn}A}
is the Weingarten operator of the rescaled surfaces
\eq{\ti M_{t}=\ti x(t,\bbS^{n}).}
For technical reasons we only want to work with normal velocities, so we introduce a time-dependent family $y(t,\cdot)\in C^{\8}(\bbS^{n},\bbS^{n})$ of diffeomorphisms in order to kill the tangent part in \eqref{Resc}. We calculate
\eq{\fr{d}{dt}\ti x(t,y(t,\cdot))=\fr{1}{F(\ti A)}\ti\nu - \fr 1n \ip{\ti x}{\ti \nu}\ti \nu-\fr 1n\ip{\ti x}{\ti \n_{j} \ti x}\ti \n_{i}\ti x \ti g^{ij} +\ti\n_{i}\ti x \dot{y}^{i}.}
Thus, if we solve the ODE system
\eq{\dot{y}^{i}=\fr 1n\ip{\ti x}{\ti\n_{j}\ti x}\ti g^{ij},}
we see that $z(t)=\ti x(t,y(t,\cdot))$ solves
\eq{\label{Reparam}\dot{z}=\br{\fr{1}{F(\ti A)}-\fr 1n \ti u}\ti \nu,}
where
\eq{\ti u=\ip{z}{\ti \nu}}
is positive due to the starshapedness of $\ti M_{t}$. This formal discussion justifies that we as well may focus on the long-time existence and regularity for the flow \eqref{Reparam}. In order to facilitate notation, we will switch back to a more convenient notation and prove the following theorem, from which \Cref{thm:ICF} then follows.

\begin{thm}\label{thm:Reparam}
Let $x_{0}$ and $F$ satisfy the assumption of \Cref{thm:ICF}. Then there exists a unique solution $x\in C^{\8}([0,\8)\x \bbS^{n},\bbR^{n+1})$ of
\eq{\label{APF}\dot{x}&=\br{\fr{1}{F(A)}-\fr un  }\nu\\
		x(0,\cdot)&=x_{0}.}
The embeddings $x(t,\cdot)$ converge smoothly to the embedding of a round sphere.
\end{thm}

\subsection*{Short time existence}
To prove that the system \eqref{APF} has a unique solution at least for a short time, we reduce it to a scalar parabolic equation and a system of ODEs. As we assume the initial hypersurface to be graphical over $\bbS^{n}$, if we already had a smooth solution for a while, the radial function would satisfy
\eq{\label{dot-rho}\dot{\rho}=\fr{\ip{x}{\dot{x}}}{\abs{x}}=\br{\fr{1}{F}-\fr{u}{n}}\fr{u}{\rho},}
as can be seen by differentiation of $\rho=\abs{x}$. From \eqref{graph-h}, \cite[Equ.~(2.4.21)]{Gerhardt:/2006} and \cite[Lemma~2.7.6]{Gerhardt:/2006} we see that $\rho=\rho(t,x^{i})$ would be the solution to the fully nonlinear equation
\eq{\label{Sc}\del_{t}{\rho}&=G(\bar\n^{2}\rho,\bar\n \rho,\rho,\cdot)\\
			\rho(0,\cdot)&=\rho_{0},}
where $\rho_{0}$ is the radial function of the initial surface $M_{0}=x(0,\bbS^{n})$ and $\bar\n$ is the Levi-Civita connection of the round metric $\s$ on $\bbS^{n}$. Also note that here $(x^{i})$ are the spherical coordinates of $x$ in the polar coordinate system of the punctured Euclidean space. The idea is to solve this Cauchy-problem, which then determines the radial functions $\rho=\rho(t,x^{i})$ of the flow hypersurfaces. Then we solve the following ODE initial value problem on $\bbS^{n}$:
\eq{\dot{x}^{i}&=\br{\fr{1}{F(A)}-\fr un  }\nu^{i}\\
		x^{i}(0)&=x_{0}^{i},}
where we note that the right hand side is fully determined by the function $\rho$ and its derivatives, which itself solely depend on $(x^{i})$. Then we plug everything together and define
\eq{x(t,\xi)=(\rho(t,x^{i}(t,\xi)),x^{i}(t,\xi)),}
which solves \eqref{APF}. In particular we note that {\it{the maximal time of existence for \eqref{APF} is entirely determined by the maximal time of existence for \eqref{Sc}}}.

It would miss the aim of this course to provide the rigorous argument behind this approach. The proof of existence for \eqref{Sc} uses solvability of linear parabolic equations in Hölder spaces and the implicit function theorem. In particular the maximal time of existence is controlled from below by estimates on the initial data. See \cite[Sec.~2.5]{Gerhardt:/2006} for some more details. We have:

\begin{thm}\label{STE}
There exists $T^{*}\leq \8$ and a unique maximal solution
\eq{x\in C^{\8}([0,T^{*})\x \bbS^{n},\bbR^{n+1})} to \eqref{APF}. If $T^{*}<\8$, then at $T^{*}$ some derivative of $x$ must blow up.
\end{thm}

\subsection*{Evolution equations}
In order to prove the immortality of the maximal solution to \eqref{APF}, by \Cref{STE} it suffices to prove uniform estimates on all derivatives of $x$. As those are controlled by derivatives of $\rho$, everything is reduced to prove regularity estimates for $\rho$.

The proof of these proceed by establishing estimates up to $C^{2}$-level as well as a lower $F$-bound by maximum principle, followed by regularity estimates for fully nonlinear parabolic operators due to Krylov and Safonov, as well as a bootstrapping argument using Schauder theory.
We need further evolution equations, which are specifically adapted to the flow \eqref{APF}. We define the operator
\eq{\cL= \del_{t}-\fr{1}{F^{2}}\tr(F'(A)\circ (\n^{2})^{\sharp})-\fr{1}{n}\ip{\rho\del_{r}}{\n^{(\cdot)}}.}

\begin{lemma}
Along the flow \eqref{APF} the
radial function $\rho=\rho(t,\xi)$ satisfies
\eq{\label{Ev-rho}\cL\rho=\fr{2}{F}\fr{u}{\rho}-\fr{\rho}{n}-\fr{1}{\rho F^{2}}\tr F'(A)+\fr{1}{\rho F^{2}}\tr(F'\circ \n\rho\otimes(\n\rho)^{\sharp}),}
while the
support function $u$ satisfies
\eq{\label{Ev-u}\cL u&=\fr{1}{F^{2}}\br{\tr(F'(A)\circ A^{2})-\fr{F^{2}}{n}}u.}
\end{lemma}

\pf{
(i)~Use \eqref{graph-h} to deduce
\eq{\tr(F'\circ (\n^{2}\rho)^{\sharp})=\fr{1}{\rho}\tr F'-\fr{u}{\rho}F-\fr{1}{\rho}\tr(F'\circ \n\rho\otimes(\n\rho)^{\sharp})}
and hence, also using \eqref{dot-rho},
\eq{\cL\rho&=\br{\fr{2}{F}-\fr{u}{n}}\fr{u}{\rho}-\fr{\rho}{n}\ip{\del_{r}}{\n^{(\cdot)}\rho}-\fr{1}{\rho F^{2}}\tr F'+\fr{1}{\rho F^{2}}\tr(F'\circ \n\rho\otimes(\n\rho)^{\sharp}).}
There holds
\eq{-\fr{1}{n}\fr{u^{2}}{\rho}-\fr{\rho}{n}\ip{\del_{r}}{\n^{(\cdot)}\rho}&=-\fr{\rho}{n} \ip{\del_{r}}{\nu}^{2}-\fr{\rho}{n}\ip{\del_{r}}{x_{\ast}\n^{(\cdot)}\rho}\\
				&=-\fr{\rho}{n}\br{ \ip{\del_{r}}{\nu}^{2}+\abs{\n\rho}^{2}}\\
				&=-\fr{\rho}{n}\br{\ip{\del_{r}}{\nu}^{2}+\sum_{i=1}^{n}\ip{\del_{r}}{\n_{i}x}^{2}}\\
				&=-\fr{\rho}{n},}
if coordinates are chosen such that $(\nu,\n_{i}x)$ is an orthonormal basis.

(ii)~The position field $r\del_{r}$ is a conformal vector field, hence for all vector fields $\bar X$ on $\bbR^{n+1}$ we have
\eq{D_{\bar X}(r\del_{r})=\bar X.}
Hence, for vector fields $X$ on $M$,
\eq{\label{Ev-u-1}\dot{u}=\ip{\dot{x}}{\nu}+\ip{\rho\del_{r}}{D_{\dot{x}}\nu}=\fr{1}{F}-\fr{u}{n}+\ip{\rho\del_{r}}{\fr{\n F}{F^{2}}}+\ip{\rho\del_{r}}{\fr{\n u}{n}},}
\eq{Xu=\ip{\rho\del_{r}}{A(X)}}
and
\eq{\label{Hess-u}\n^{2}u(X,Y)&=Y(Xu)-(\n_{Y}X)u=h(X,Y)-h(X,A(Y))u+\ip{\rho\del_{r}}{\n_{Y}A(X)}.}
The result follows from combining these equalities, also using the Codazzi equation to cancel the $\n F$-terms and the homogeneity of $F$ which implies
\eq{\tr(F'(A)\circ A)=F.}
}

We also need specific curvature evolution equations to estimate the principal curvatures and $F$ from below.

\begin{lemma}\label{W-par}
The Weingarten operator satisfies,
\eq{\label{Ev-W-2}\cL A=\fr{1}{F^{2}}(F'\circ A^{2})A-\fr{2A^{2}}{F}+\fr{A}{n}-\fr{2}{F^{3}}\n F\otimes (\n F)^{\sharp}+\fr{1}{F^{2}}d^{2}F(\n_{(\cdot)}A,\n^{(\cdot)}A),}
while the curvature function $F$ satisfies
\eq{\cL F&=-\fr{1}{F^{2}}\br{\tr(F'(A)\circ A^{2})-\fr{F^{2}}{n}}F-\fr{2}{F^{3}}(F'\circ \n F\otimes (\n F)^{\sharp}).}
\end{lemma}

\pf{
(i)~From \eqref{Ev-W-1} we calculate
\eq{\label{W-par-1}\dot{A}&=\n\grad \br{\fr{u}{n}-\fr{1}{F}}+ \br{\fr{u}{n}-\fr{1}{F}}				A^{2}\\
			&=\fr{(\n^{2} u)^{\sharp}}{n}+\fr{(\n^{2} F)^{\sharp}}{F^{2}}-\fr{2}{F^{3}}\n F\otimes (\n F)^{\sharp}+ \br{\fr{u}{n}-\fr{1}{F}}A^{2}\\
			&=\fr{A}{n}+\fr{1}{n}\ip{\rho\del_{r}}{\n^{(\cdot)}A}+\fr{(\n^{2} F)^{\sharp}}{F^{2}}-\fr{2}{F^{3}}\n F\otimes (\n F)^{\sharp}-\fr{1}{F}A^{2}.
		}
We have to analyze the term $\n^{2}F$ and do this is a local coordinate frame. There hold
\eq{\n_{i}F=dF(A)\n_{i}A}
and
\eq{\n_{ji}F=d^{2}F(A)(\n_{i}A,\n_{j}A)+dF(A)\n_{ji}A.}
We have to swap indices in $\n_{ij}A=\n_{ij}h^{k}_{l}$.

\eq{
\n_{ji}h_{l}^{k}&=\n_{jl}h_{i}^{k}\\
                &=\n_{lj}h_{i}^{k}+{R_{jla}}^kh_{i}^{a}-{R_{jli}}^a h_{ka}\\
                &=\n^{k}_{l}h_{ij}+{R_{jla}}^kh_{i}^{a}-{R_{jli}}^a h_{ka}\\
                &=\n^{k}_{l}h_{ij}+(h^{k}_{j}h_{la}-h_{ja}h^{k}_{l})h^{a}_{i}-(h^{a}_{j}h_{li}-h^{a}_{l}h_{ij})h_{ka}.
             }
Applying $d F=d F(A)=(F^{l}_{k})$ to this, while using the $1$-homogeneity and that $dF(A)$ commutes with $A$, gives
\eq{F^{l}_{k}\n_{ji}h^{k}_{l}&=F^{l}_{k}\n^{k}_{l}h_{ij}+F^{l}_{k}(h^{k}_{j}h_{la}-h_{ja}h^{k}_{l})h^{a}_{i}-F^{l}_{k}(h^{a}_{j}h_{li}-h^{a}_{l}h_{ij})h_{ka}\\
					&=F^{l}_{k}\n^{k}_{l}h_{ij}-F^{l}_{k}h_{ja}h^{k}_{l}h^{a}_{i}+F^{l}_{k}h^{a}_{l}h_{ij}h_{ka}\\
					&=F^{l}_{k}\n^{k}_{l}h_{ij}-Fh_{ja}h^{a}_{i}+F^{l}_{k}h^{a}_{l}h_{ka}h_{ij}.
					}
Application of the sharp-operator gives
\eq{(\n^{2}F)^{\sharp}&=d^{2}F(\n_{(\cdot)}A,\n^{(\cdot)}A)+\tr(F'\circ (\n^{2}A)^{\sharp})-FA^{2}+(F'\circ A^{2})A.				
}
Inserting this into \eqref{W-par-1} gives the first equation.

(ii)~To get the equation for $F$ calculate in local coordinates
\eq{\dot{F}=F^{l}_{k}\dot{h}^{k}_{l}}
and use
\eq{F^{l}_{k}F^{i}_{j}\n^{j}_{i}h^{k}_{l}=F^{i}_{j}\n^{j}_{i}F-F^{i}_{j}d^{2}F(A)(\n_{i}A,\n^{j}A).}
}

\subsection*{A priori estimates}

The following estimates control the flow up to $C^{2}$-level for the function
\eq{\rho\cn [0,T^{*})\x\bbS^{n}\ra \bbR.} The following proof contains some common of tricks on how to estimate solutions to parabolic equations. It should be interesting even outside the world of curvature flows.

\begin{lemma}\label{AP}
There exists a constant $c>0$, which only depends on the initial hypersurface, such that
\enum{
\item \eq{\label{AP-1}\min_{\bbS^{n}}\rho(0,\cdot)\leq \rho\leq \max_{\bbS^{n}}\rho(0,\cdot),}
\item \eq{\label{AP-2}c^{-1}\leq u\leq c,}
\item \eq{\label{AP-3}c^{-1}\leq F\leq c,}
\item \eq{\label{AP-4}\abs{A}^{2}\leq c.}
}
It follows that there exists a compact set $K\sub \G$, in which the principal curvatures range during the whole evolution.
\end{lemma}

\pf{
(i)~Define
\eq{\ti\rho(t)=\max_{\bbS^{n}}\rho(t,\cdot).}
Then $\ti\rho$ is Lipschitz and hence differentiable almost everywhere. It can be shown that at points of differentiability there holds
\eq{\fr{d}{dt}{\ti\rho}=\dot\rho(t,\xi_{t}),}
where $\xi_{t}$ is a point where the maximum is attained. This technical argument is due to Hamilton \cite{Hamilton:/1986}. From \eqref{dot-rho} we get (note $d\rho=0$)
\eq{\dot{\ti\rho}=\fr{1}{F}-\fr{\ti\rho}{n}.}
Now we recall that $F$ depends on the second fundamental form which is related to $\rho$ via \eqref{ThetaDeriv}, which gives \eq{A=\fr{\id}{\rho}-(\n^{2}\rho)^{\sharp}\geq \fr{\id}{\rho}}
at $\xi_{t}$. Hence at $\xi_{t}$ we have
\eq{F(A)\geq \fr{n}{\ti\rho},}
and hence ${\ti\rho}$ is non-increasing. The same argument at minimal points gives that the minimum of $\rho$ is non-decreasing, which concludes the argument.

(ii)~In order to bound $u^{-1}$, we first note that $uF$ is bounded from above and below, which can be seen as follows. Define
\eq{w=\log u+\log F,}
then $w$ satisfies
\eq{\cL w=\fr{1}{u^{2}F^{2}}\tr(F'\circ \n u\otimes (\n u)^{\sharp})-\fr{1}{F^{4}}\tr(F'\circ \n F\otimes (\n F)^{\sharp}).}
At critical points of $w$ there holds
\eq{\fr{\n u}{u}=-\fr{\n F}{F}.}
Hence, as above, the functions $\max w$ and $\min w$ are non-increasing/decreasing. We can use the boundedness of $uF$ to prove that $u^{-1}$ is bounded as well. This and the subsequent estimates all boil down to finding appropriate test functions.

The evolution equation of $u^{-1}$ has one bad positive term, which prevents us from estimating it directly. Namely there holds
\eq{\cL u^{-1}\leq \fr{u^{-1}}{n}-\fr{2}{u^{3}F^{2}}\tr(F'\circ \n u\otimes(\n u)^{\sharp}).}
However, we already have one bounded quantity, $\rho$, and we can use it to build test functions. Define
\eq{w=\log u^{-1}+\la \rho,\q \la>0.}
There holds, due to $Fu\geq c>0$,
\eq{\cL w\leq \fr{1}{n}  +\la\fr{c}{\rho}u^{2}-\fr{\la\rho}{n}<0}
at all critical points of $w$ where $w$ is large enough, provided $\la$ is chosen large enough. Hence $w$ is bounded. In turn $u^{-1}$ is bounded. The upper bound for $u$ simply follows from
\eq{u\leq \rho.}

(iii)~Follows directly from the bounds on $u$ and those on $uF$.

(iv)~We use
\eq{\dot{g}=2\br{\fr{1}{F}-\fr{u}{n}} h}
to deduce
\eq{\cL h=\cL h^{k}_{i}g_{kj}+2\br{\fr{1}{F}-\fr un}h_{kj}h^{k}_{i} }
and hence
\eq{\cL h=\fr{F'\circ A^{2}}{F^{2}}h-\fr{2u}{n}h(A,\cdot)+\fr{h}{n}-\fr{2}{F^{3}}\n F\otimes \n F+\fr{1}{F^{2}}d^{2}F(\n_{(\cdot)}A,\n_{(\cdot)}A).}
The only angry looking term in the evolution of $h$ is the first one. As it also appears in the evolution of $u$, we cancel it with this one.
Suppose the function
\eq{z=u^{-1}\ka_{n}}
attains a maximal value at a point $(t_{0},\xi_{0})$.
Let $\eta\in T_{\xi_{0}}\bbS^{n}$ be an eigenvector corresponding to $\ka_{n}$ and extend $\eta$ locally to a vector field such that $\n \eta(t_{0},\xi_{0})=0$.
Define
\eq{w=\fr{h(\eta,\eta)}{g(\eta,\eta)}u^{-1}.}
Then, locally around $(t_{0},\xi_{0})$ there holds 
\eq{w\leq z,\q w(t_{0},\xi_{0})=z(t_{0},\xi_{0}).}
Hence $w$ also attains a local maximum at this point and it suffices to locally estimate $w$.
At $(t_{0},\xi_{0})$ there holds
\eq{\cL w&\leq-\fr{2}{n}\fr{h(A(\eta),\eta)}{g(\eta,\eta)}+\fr{2}{n}w-2\fr{h(\eta,\eta)^{2}}{g(\eta,\eta)^{2}}\br{\fr{1}{Fu}-\fr{1}{n}}\\
		&=\fr{2}{n}w-\fr{2u}{F}w^{2}-\fr{2}{n}\fr{h(A(\eta),\eta)g(\eta,\eta)-h(\eta,\eta)^{2}}{g(\eta,\eta)^{2}}\\
		&=\fr{2}{n}w-\fr{2u}{F}w^{2}-\fr{2}{n}\fr{\abs{A(\eta)}^{2}\abs{\eta}^{2}-g(A(\eta),\eta)^{2}}{\abs{\eta}^{2}},}
which is negative for large $w$ due to Bunjakowski-Cauchy-Schwarz and where we used the concavity of $F$. Hence $w$ is bounded and thus all eigenvalues of $A$ are bounded from above. As $F$ is also bounded from below, we deduce from the concavity of $F$ that
\eq{0<F\leq H,}
\cite[Lemma~2.2.20]{Gerhardt:/2006}. Hence
\eq{\ka_{1}>(1-n)\ka_{n}\geq -c}
and we obtain $\abs{A}^{2}\leq c$. If there existed a sequence $\ka(t_{n},\xi_{n})\in \G$ that leaves every compact set of $\G$, any subsequential limit of this sequence would lie on $\del\G$, which is impossible due to
\eq{F(t_{n},\xi_{n})\geq c^{-1},\q F_{|\del\G}=0.}
}

As a corollary we obtain full spatial $C^{2}$-estimates for the radial function $\rho$.

\begin{cor}
There exists a constant $c$, which only depends on the initial hypersurface, such that
\eq{\abs{\rho(t,\cdot)}_{C^{2}(\bbS^{n})}\leq c\q\fa t\in [0,T^{*}).}
\end{cor}

\pf{
The $C^{0}$-bound of $\rho$ follows from \eqref{AP-1}. As we are dealing with graphs over $\bbS^{n}$ in the product space
\eq{\bbR^{n+1}\bs\{0\}=(0,\8)\x\bbS^{n}, \q \ip{\cdot}{\cdot}=dr^{2}+r^{2}\s,}
the normal $\nu(\rho(t,\xi))$ is given by
\eq{\nu=\fr{(1,-\rho^{-2}\s^{ik}\del_{i}\rho)}{\rt{1+\rho^{-2}\abs{d\rho}^{2}_{\s}}}}
and hence the support function is
\eq{u=\rho\ip{\del_{r}}{\nu}=\fr{\rho}{\rt{1+\rho^{-2}\abs{d\rho}^{2}_{\s}}}.}
As $\rho$ and $u$ are uniformly bounded, so is $\abs{d\rho}_{\s}$, which gives the $C^{1}$-estimate. $C^{2}$-estimates follow from curvature estimates and the representation of the second fundamental form in terms of the second derivatives of $\rho$, \eqref{ThetaDeriv}.
}

The key for higher order estimates is a regularity result due to Krylov \cite{Krylov:/1987}. We state a very accessible formulation of this result as it can be found in a note by Ben Andrews \cite[Thm.~4]{Andrews:/2004}.

\begin{thm}\label{Kr}
Let $\Om\sub\bbR^{n}$ be open and suppose $\rho\in C^{4}((0,T]\x\Om)$ satisfies
\eq{\del_{t}\rho=G(D^{2}\rho,D\rho,\rho,\cdot),}
where $G$ is concave in the first variable. Then for any $\tau>0$ and $\Om'\Subset \Om$ there holds
\eq{&\sup_{s,t\in [\tau,T], p,q\in \Om'}\br{\fr{\abs{D^{2}\rho(p,t)-D^{2}\rho(q,t)}}{\abs{p-q}^{\al}+\abs{s-t}^{\fr{\al}{2}}}+\fr{\abs{\del_{t}\rho(p,t)-\del_{t}\rho(q,t)}}{\abs{p-q}^{\al}+\abs{s-t}^{\fr{\al}{2}}}}\\
		\hp{=}+~&\sup_{s,t\in[\tau,T],p\in\Om'}\fr{\abs{D\rho(p,t)-D\rho(p,s)}}{\abs{s-t}^{\fr{(1+\al)}{2}}}\leq C,}
where $\al$ depends on $n$ and the ellipticity constants $\la,\La$ of $F'$, and $C$ depends on $n$, $\la,\La$, bounds for $\abs{D^{2}\rho}$ and $\abs{\del_{t}\rho}$, $d(\Om',\del\Om)$, $\tau$ and the bounds on the other first and second derivatives of $G$.
\end{thm}

 As $\rho$ satisfies the fully nonlinear equation
\eq{\del_{t}\rho=G(\bar\n^{2}\rho,\bar\n\rho,\rho,\cdot)=\br{\fr{1}{F(A)}-\fr{u}{n}}\fr{\rho}{u},}
cf. \cite[Equ.~(2.4.21)]{Gerhardt:/2006},
let us quickly check the assumptions of this theorem are satisfied. We use
\eqref{graph-h}, \cite[Lemma~2.7.6]{Gerhardt:/2006} and \cref{AP} to obtain the uniform ellipticity of
\eq{\fr{\del G}{\del \rho_{ij}}=-\fr{1}{F^{2}}dF \fr{\del A}{\del \rho_{ij}}}
and the convexity of $G$ in the first variable,
\eq{\fr{\del^{2}G}{\del\rho_{ij}\del\rho_{kl}}=\fr{2}{F^{3}}dF\fr{\del A}{\del \rho_{ij}}dF \fr{\del A}{\del \rho_{kl}}-\fr{1}{F^{2}}d^{2}F \br{\fr{\del A}{\del \rho_{ij}}, \fr{\del A}{\del \rho_{kl}}}}
where we used that $A$ depends on $\bar\n^{2}\rho$ linearly. Hence \cref{Kr} does not apply directly, but we see that $-\rho$ satisfies an equation with a concave operator, to which we can apply the theorem.
Hence $\rho$ lies in the parabolic Hölder space $H^{2+\al,\tfrac{2+\al}{2}}([0,T]\x\bbS^{n})$ for every $T<T^{*}$ with estimates independent of $T$. A standard bootstrapping argument using parabolic Schauder estimates implies uniform $C^{k}$-estimates of $\rho$ for every $k$. It follows:

\begin{cor}
There exists a constant $c$, depending only on initial data and $k$, such that
\eq{\abs{x(t,\cdot)}_{C^{k}(\bbS^{n})}\leq c \q\fa 0\leq t<T^{*}.}
The solution to \eqref{APF} is immortal.
\end{cor}

\pf{We have already seen the argument for the uniform estimates. The argument for immortality of the solution goes as follows. Suppose $T^{*}<\8$. From \Cref{STE} we know that the maximal time of existence can be estimated from below in terms of estimates for the initial data. As we have uniform estimates up to $T^{*}$, we may move as close to $T^{*}$ as required to exceed $T^{*}$ once we start with
\eq{\ti M_{0}=M_{t},\q T^{*}-\e<t<T^{*},}
where $\e$ is chosen such that the flow with initial data $\ti M_{0}$ exists longer than $\e$. Due to uniqueness we can extend our original flow and thus have shown that $T^{*}=\8$.
}

\subsection*{Convergence to a round sphere}
To conclude the proof of \cref{thm:Reparam}, we have to show convergence of the embeddings $x(t,\cdot)$ to the embedding of a round sphere. We use the strong maximum principle.
\begin{thm}
The solution $x$ to \eqref{APF} limits to the embedding of a round sphere as $t\ra \8$.
\end{thm}

\pf{
We have shown that the radial function $\rho$ satisfies a uniformly parabolic equation. Hence its oscillation
\eq{\osc \rho(t)=\max_{\bbS^{n}}\rho(t,x(t,\cdot))-\min_{\bbS^{n}}\rho(t,x(t,\cdot))=\max_{\bbS^{n}}x^{0}(t,\cdot)-\min_{\bbS^{n}}x^{0}(t,\cdot)}
 is strictly decreasing, unless it is zero. Suppose it would not converge to zero as $t\ra \8$. Then it converges to some other value
 \eq{\osc\rho(t)\ra c_{0}>0,\q t\ra \8. }
Due to our uniform estimates, a diagonal argument and Arzela-Ascoli, the sequence of flows
\eq{x_{k}(t,\xi):=x(t+k,\xi)}
subsequentially converges to a limit flow $x_{\8}$ with corresponding radial function $\rho_{\8}$. There holds
\eq{\osc\rho_{\8}(t)=\lim_{k\ra \8}\osc \rho_{k}(t)=c_{0}.  }
Hence the strong maximum principle holding for $\rho_{\8}$ is violated if $c_{0}>0$. Thus $\osc\rho(t)\ra 0$ as $t\ra 0$ and hence the flow converges to a sphere centered at the origin.
}

\section{Alexandrov-Fenchel inequalities}
We use \Cref{thm:Reparam} to prove the classical Alexandrov-Fenchel inequalities for starshaped hypersurfaces with $\s_{k}>0$, cf. \cite{GuanLi:08/2009}. These are inequalities between so-called {\it{higher order volumes}}. To motivate the terminology, let us consider a convex body, i.e. a compact convex $K$ set with non-empty interior and its {\it{$\e$-parallel body}}
\eq{K_{\e}=\{x\in \bbR^{n+1}\cn \dist(K,x)\leq \e\}.}
A classical result is {\it{Steiner's formula}}, which provides a Taylor expansion of the volume of $K_{\e}$:
\eq{\vol(K_{\e})=\sum_{k=0}^{n}\binom{n+1}{k}W_{k}(K)\e^{k}\q \fa \e\geq 0,}
where the $W_{k}(K)$ are called the quermassintegrals of $K$, cf. \cite{Schneider:/2014}. Locally, such an expansion even holds for non-convex domains. In the following we prove this and a useful representation formula. First we need a general variational formula, where $S_{k}$ is the operator function associated to the elementary symmetric polynomial $s_{k}$, see \eqref{ES}. We also define
\eq{s_{0}:=1,\q s_{-1}:=u.}

We will use the following facts about the $S_{k}$ without proof:

\eq{dS_{k}A=kS_{k} \q \fa 0\leq k\leq n,}
\eq{ \tr(dS_{k+1})=(n-k)S_{k}\q \fa 0\leq k\leq n-1,}
\eq{dS_{k}A^{2}=S_{1}S_{k}-(k+1)S_{k+1}\q \fa 0\leq k\leq n-1.}
Furthermore $dS_{k}$ is divergence free.
Hence we can deduce:
\begin{lemma}\label{Sk-Ev}
Let $x$ and $\cF$ as in \Cref{Ev} with $M$ compact. For every $1\leq k\leq n$ there holds
\eq{\del_{t}\int_{M_{t}}S_{k-1}=-k\int_{M_{t}}\cF S_{k}.}
For $k=0$ there holds
\eq{\del_{t}\int_{M_{t}}\ip{x}{\nu}=-(n+1)\int_{M_{t}}\cF.}
\end{lemma}

\pf{
For $k=0$ we have
\eq{\del_{t}\int_{M_{t}}\ip{x}{\nu}&=-\int_{M_{t}}\cF+\int_{M_{t}}(\ip{x}{\grad \cF}-\cF H\ip{x}{\nu})\\
						&=-\int_{M_{t}}\cF+\int_{M_{t}}\dive_{M_{t}}(\cF x^{\top})-n\int_{M_{t}}\cF.}

For $k=1$ we have
\eq{\del_{t}\Area(M_{t})=-\int_{M_{t}}\cF H=-\int_{M_{t}}\cF S_{1},}
while for $2\leq k\leq n$ we calculate using \eqref{Ev-W-1}:
\eq{\del_{t}\int_{M_{t}}S_{k-1}&=-\int_{M_{t}}S_{k-1}\cF S_{1}+\int_{M_{t}}\tr(dS_{k-1}\circ \n \grad \cF)+\int_{M_{t}}\cF \tr(dS_{k-1}\circ A^{2})\\
			&=\int_{M_{t}}\cF (dS_{k-1}A^{2}-S_{1}S_{k-1})\\
			&=-k\int_{M_{t}}\cF S_{k}.}
}

Now we can prove a local Steiner's formula for $C^{2}$-domains.

\begin{lemma}
Let $\Om\sub\bbR^{n+1}$ be a bounded domain with $C^{2}$-boundary and let $\bar\Om_{\e}$ be the $\e$-parallel body. Then there exists some $\e_{0}>0$ such that for all $0\leq \e<\e_{0}$ we have the expansion
\eq{\label{Steiner}\vol(\bar\Om_{\e})=\sum_{k=0}^{n}\binom{n+1}{k}W_{k}(\Om)\e^{k},}
where $W_{0}(\Om)=\vol(\Om)$ and 
 \eq{W_{k}(\Om)=\fr{1}{(n+1)\binom{n}{k-1}}\int_{\del\Om}s_{k-1}(\ka_{i}),\q 1\leq k\leq n+1.}
\end{lemma}

\pf{
There holds
\eq{\label{vol}\vol(\bar\Om_{\e})=\fr{1}{n+1}\int_{\Om_{\e}}\dive x=\fr{1}{n+1}\int_{\del\bar\Om_{\e}}\ip{x}{\nu_{\e}}}
and hence
\eq{W_{0}(\Om)=\vol(\Om)=\fr{1}{n+1}\int_{\del\Om}s_{-1}.}

The parallel hypersurfaces $\del\bar\Om_{\e}$, which are $C^{2}$-hypersurfaces for small $\e$, can be seen as the flow hypersurfaces of the flow
\eq{\del_{\e}x=\nu_{\e}.}
According to \Cref{Sk-Ev} we obtain
\eq{\del_{\e}\vol(\Om_{\e})=\Area(\del\Om_{\e})}
and
\eq{\del_{\e}^{k}\vol(\Om_{\e})=\del^{k-1}_{\e}\int_{\del\Om_{\e}}1=(k-1)!\int_{\del\Om_{\e}}S_{k-1}\q \fa 1\leq k\leq n+1.}
For $k> n+1$ there holds
\eq{\del_{\e}^{k}\vol(\Om_{\e})=0}
due to Gauss-Bonnet.
Defining the $W_{k}$ according to the Taylor expansion in \eqref{Steiner} we see that they must have the form
\eq{W_{k}(\Om)=\fr{1}{k!\binom{n+1}{k}}\del_{\e}^{k}\vol(\bar\Om_{\e})_{|\e=0}=\fr{1}{k\binom{n+1}{k}}\int_{\del\Om}S_{k-1},}
which is the claimed formula.
}

Hence the $W_{k}(\Om)$ are nothing but coefficients of higher order in the Taylor expansion of volume with respect to fattening of the boundary. The isoperimetric inequality provides an estimate between $W_{0}$ and $W_{1}$ and hence it is natural to ask whether such an estimate also holds between the other higher order volumes. While for convex bodies such estimates have long been known, see for example \cite{Schneider:/2014} for a broad overview, here we want to use \Cref{thm:Reparam} to prove them for starshaped hypersurface with a certain curvature condition. This approach is due to Pengfei Guan and Junfang Li \cite{GuanLi:08/2009}.

\begin{defn}
A domain $\Om\sub\bbR^{n+1}$ is called {\it{k-convex}}, if throughout $\del\Om$ the principal curvatures lie in the closure of the cone
\eq{\G_{k}=\{\ka\in \bbR^{n}\cn s_{m}(\ka)>0\q \fa m\leq k\}.}
$\Om$ is called {\it{strictly k-convex}}, if the principal curvatures lie in $\G_{k}$.
\end{defn}

\begin{thm}
Let $\Om\sub\bbR^{n+1}$ be a starshaped and $k$-convex domain, then there holds
\eq{\fr{W_{k+1}(\Om)}{W_{k+1}(B)}\geq \br{\fr{W_{k}(\Om)}{W_{k}(B)}}^{\fr{n-k}{n+1-k}},}
where $B$ is the unit ball in $\bbR^{n+1}$. Equality holds precisely if $\Om$ is a ball.
\end{thm}

\pf{
We use \cref{thm:Reparam} with
\eq{F=n\fr{\binom{n}{k-1}}{\binom{n}{k}}\fr{s_{k}}{s_{k-1}}=\fr{nk}{n-k+1}\fr{s_{k}}{s_{k-1}}}
and start the flow with $M_{0}=\del\Om$.
Due to the $k$-convexity of $\Om$ the assumptions of this theorem are satisfied. We calculate that along the flow
\eq{\dot{x}=\br{\fr{1}{F}-\fr{u}{n}}\nu}
there holds
\eq{\del_{t}W_{k}(\Om_{t})=\fr{k}{(n+1)\binom{n}{k-1}}\int_{M_{t}}\br{\fr{n-k+1}{nk}\fr{s_{k-1}}{s_{k}}-\fr{u}{n}}s_{k}.}
As the $dS_{k}$ are divergence free, we obtain after tracing \eqref{ThetaDeriv} with respect to $dS_{k}$ and integration:
\eq{(n-k+1)\int_{M_{t}}s_{k-1}=k\int_{M_{t}}us_{k}}
and thus
\eq{\del_{t}W_{k}(\Om_{t})=0.}
On the other hand we obtain
\eq{\fr{n(n+1)\binom{n}{k}}{k+1}\del_{t}W_{k+1}(\Om_{t})&=\int_{M_{t}}\br{\fr{n-k+1}{k}\fr{s_{k-1}s_{k+1}}{s_{k}}-us_{k+1}}\\
				&=\int_{M_{t}}\br{\fr{n-k+1}{k}\fr{s_{k-1}s_{k+1}}{s_{k}}-\fr{n-k}{k+1}s_{k}}\\
				&\leq 0,}
by the Newton-Maclaurin inequalities. Hence $W_{k+1}$ is decreasing along the flow. We know the flow converges to a sphere, on which the desired inequality holds with equality. Hence on $\Om$ the inequality is valid. The equality case follows, since the Newton-Maclaurin inequalities hold with equality precisely in umbilical points. This concludes the proof.
}

\bibliographystyle{shamsplain}
\bibliography{Bibliography}
\end{document}